\documentclass[11pt]{amsart}
\usepackage{amsmath,amsthm,amsfonts,amssymb,latexsym,epsfig}

\newtheorem{theorem}{Theorem}[section]

\newtheorem{lemma}{Lemma}[section]

\numberwithin{equation}{section}

\theoremstyle{definition}

\theoremstyle{remark}

\begin{document}
\title{A Note on Carleman's Inequality}
\author{Peng Gao}
\address{Department of Computer and Mathematical Sciences,
University of Toronto at Scarborough, 1265 Military Trail, Toronto
Ontario, Canada M1C 1A4} \email{penggao@utsc.utoronto.ca}
\date{June 15, 2007.}


\begin{abstract}
 We study a weighted version of Carleman's inequality via Carleman's original approach. As an application of our result, we prove a conjecture of Bennett.
\end{abstract}

\maketitle
\section{Introduction}
\label{sec 1} \setcounter{equation}{0}
  
   The well-known Carleman's inequality asserts that for convergent infinite series $\sum a_n$ with non-negative terms, one has
\begin{equation*}
   \sum^\infty_{n=1}(\prod^n_{k=1}a_k)^{\frac 1{n}}
\leq e\sum^\infty_{n=1}a_n,
\end{equation*}
   with the constant $e$ best possible.

  There is a rich literature on many different proofs of Carleman's inequality as well as its generalizations and extensions. We shall refer the readers to the survey articles  \cite{P&S} and \cite{D&M} as well as the references therein for
an account of Carleman's inequality. 

  From now on we will assume $a_n \geq 0$ for $n \geq 1$ and any
   infinite sum converges. Our goal in this paper is to study the following weighted Carleman's inequality:
\begin{equation}
\label{1}
   \sum^\infty_{n=1}G_n
\leq C\sum^\infty_{n=1}a_n,
\end{equation}
  where
\begin{equation}
\label{2}
   G_n=\prod^n_{k=1}a^{\lambda_k/\Lambda_n}_k, \hspace{0.1in} \Lambda_n=\sum^n_{k=1}\lambda_k, ~~\lambda_k \geq 0, ~~\lambda_1>0.
\end{equation}
   The task here is to determine the best constant $C$ so that inequality \eqref{1} holds for any non-negative sequence $\{a_n \}^{\infty}_{n=1}$.
  
   One approach to our problem here is to deduce inequality \eqref{1} via  $l^{p}$ operator norm of the corresponding weighted mean matrix. We recall here that a matrix $A=(a_{j,k})$ is said to be a weighted mean matrix if its entries satisfy:
\begin{equation}
\label{3}
    a_{j,k}=\lambda_k/\Lambda_j, 1 \leq k \leq j; \hspace{0.1in} a_{j,k}=0, k>j,
\end{equation}
   where the notations are as in \eqref{2}.
   For $p>1$, let $l^p$ be the Banach space of all complex sequences ${\bf b}=(b_n)_{n \geq 1}$ with norm
\begin{equation*}
   ||{\bf b}||: =(\sum_{n=1}^{\infty}|b_n|^p)^{1/p} < \infty.
\end{equation*}
   The $l^{p}$ operator norm $||A||_{p,p}$ of $A$ for $A$ as defined in \eqref{3} is then defined as the $p$-th root of the smallest value of the
   constant $U$ so that the following inequality holds for any ${\bf b} \in l^p$:
\begin{equation}
\label{4}
   \sum^{\infty}_{n=1}\big{|}\sum^{\infty}_{k=1}\lambda_kb_k/\Lambda_n
   \big{|}^p \leq U \sum^{\infty}_{n=1}|b_n|^p.
\end{equation}
   
   In an unpublished dissertation \cite{Car}, Cartlidge studied
weighted mean matrices as operators on $l^p$ and obtained the
following result (see also \cite[p. 416, Theorem C]{B1}).
\begin{theorem}
\label{thm02}
    Let $1<p<\infty$ be fixed. Let $A=(a_{j,k})$ be a weighted mean matrix given by \eqref{3}. If
\begin{equation}
\label{022}
    L=\sup_n(\frac {\Lambda_{n+1}}{\lambda_{n+1}}-\frac
    {\Lambda_n}{\lambda_n}) < p ~~,
\end{equation}
    then
    $||A||_{p,p} \leq p/(p-L)$.
\end{theorem}

   The above theorem implies that one can take $U=(p/(p-L))^p$ in inequality \eqref{4} for any weighted mean matrix $A$ satisfying \eqref{022}. We note here by a change of variables $b_k \rightarrow a^{1/p}_k$ in \eqref{4} and on letting $p \rightarrow +\infty$, one obtains inequality \eqref{1} with $C=e^{L}$ as long as \eqref{022} is satisfied with $p$ replaced by $+\infty$ there.
   
   In this note, we will study inequality \eqref{1} via Carleman's original approach and we shall prove in the next section the following:
\begin{theorem}
\label{thm1}
  Suppose that
\begin{equation}
\label{5}
  M=\sup_n\frac
    {\Lambda_n}{\lambda_n}\log \Big(\frac {\Lambda_{n+1}/\lambda_{n+1}}{\Lambda_n/\lambda_n} \Big ) < +\infty,
\end{equation}
  then inequality \eqref{1} holds with $C=e^M$.
\end{theorem}
 
   We point out here that the result of Theorem \ref{thm1} is better than what one can deduce from Cartlidge's result as discussed above. This can be seen by noting that \eqref{5} is equivalent to 
\begin{equation*}
  \frac {\Lambda_{n+1}\lambda_{n}}{\Lambda_{n}\lambda_{n+1}} \leq e^{\lambda_nM/\Lambda_n},
\end{equation*}
  for any integer $n \geq 1$. Suppose now \eqref{022} is satisfied, then the case $n=1$ of \eqref{022} implies $L>0$ and it is easy to check that
\begin{equation*}
   \frac {\Lambda_{n+1}\lambda_{n}}{\Lambda_{n}\lambda_{n+1}} = 1+\frac
    {\lambda_n}{\Lambda_n}(\frac {\Lambda_{n+1}}{\lambda_{n+1}}-\frac
    {\Lambda_n}{\lambda_n}) \leq 1+\frac
    {\lambda_n}{\Lambda_n}L \leq e^{\lambda_nL/\Lambda_n},
\end{equation*} 
  from which we deduce that $M \leq L$.

  Bennett \cite[p. 829]{Be1} conjectured that inequality \eqref{1} holds for $\lambda_k=k^{\alpha}$ for $\alpha > -1$ with $C=1/(\alpha+1)$. As the cases $-1 < \alpha \leq 0$ or $\alpha \geq 1$ follow directly from Cartlidge's result above (Theorem \ref{thm02}), the only case left unknown is when $0< \alpha <1$. As an application of Theorem \ref{thm1}, we shall prove Bennett's conjecture in Section \ref{sec 3}.
   
\section{Proof of Theorem \ref{thm1}}
\label{sec 2} \setcounter{equation}{0}
  It suffices to establish our assertion with the infinite summation in \eqref{1} replaced by any finite summation, say from $1$ to $N \geq 1$ here. We now follow Carleman's approach by determing the maximamum value $\mu_N$ of $\sum^N_{n=1}G_n$ subject to the constraint $\sum^N_{n=1}a_n=1$ using Lagrange multipliers. It is easy to see that we may assume $a_n > 0$ for all $1 \leq n \leq N$ when the maximamum is reached. We now define
\begin{equation*}
  F({\bf a}; \mu)=\sum^N_{n=1}G_n-\mu (\sum^N_{n=1}a_n-1),
\end{equation*}
  where ${\bf a}=(a_n)_{1 \leq n \leq N}$. By the Lagrange method, we have to solve $\nabla F=0$, or the following system of equations:
\begin{equation}
\label{2.1}
  \mu a_k=\sum^N_{n=k}\frac {\lambda_kG_n}{\Lambda_n}, \hspace{0.1in} 1 \leq k \leq N; \hspace{0.1in} \sum^N_{n=1}a_n=1.
\end{equation}
   We note that on summing over $1 \leq k \leq N$ of the first $N$ equations above, we get
\begin{equation*}
  \sum^N_{n=1}G_n=\mu.
\end{equation*}
   Hence we have $\mu=\mu_N$ in this case which allows us to recast the equations \eqref{2.1} as:
\begin{equation*}
  \mu_N \frac {a_k}{\lambda_k}=\sum^N_{n=k}\frac {G_n}{\Lambda_n}, \hspace{0.1in} 1 \leq k \leq N; \hspace{0.1in} \sum^N_{n=1}a_n=1.
\end{equation*}
  On subtracting consecutive equations, we can rewrite the above system of equations as:
\begin{equation*}
  \mu_N (\frac {a_k}{\lambda_k}-\frac {a_{k+1}}{\lambda_{k+1}})=\frac {G_k}{\Lambda_k}, \hspace{0.1in} 1 \leq k \leq N-1; \hspace{0.1in}  \mu_N \frac {a_N}{\lambda_N}=\frac {G_N}{\Lambda_N}; \hspace{0.1in} \sum^N_{n=1}a_n=1.
\end{equation*}
   
   Now we define for $1 \leq k \leq N-1$,
\begin{equation*}
  \omega_k = \frac {\Lambda_k}{\lambda_k}-\frac {\Lambda_k a_{k+1}}{\lambda_{k+1}a_k},
\end{equation*}
  so that we can further rewrite our system of equations as:
\begin{equation*}
  \mu_N a_k \omega_k=G_k, \hspace{0.1in} 1 \leq k \leq N-1; \hspace{0.1in}  \mu_N \frac {a_N}{\lambda_N}=\frac {G_N}{\Lambda_N}; \hspace{0.1in} \sum^N_{n=1}a_n=1.
\end{equation*}
   It is easy to check that for $1 \leq k \leq N-2$,
\begin{equation*}
   \omega^{\Lambda_{k+1}}_{k+1}=\frac 1{ \mu^{\lambda_{k+1}}_N}\Big (\frac {\omega_{k}}{\frac {\lambda_{k+1}}{\Lambda_k}(\Lambda_k/\lambda_k-\omega_{k})} \Big )^{\Lambda_{k}}.
\end{equation*}
   We now define a sequence of real functions $\Omega_k(\mu)$ inductively by setting $\Omega_1(\mu)=1/\mu$ and
\begin{equation}
\label{2.2}
   \Omega^{\Lambda_{k+1}}_{k+1}(\mu)=\frac 1{ \mu^{\lambda_{k+1}}}\Big (\frac {\Omega_{k}}{\frac {\lambda_{k+1}}{\Lambda_k}(\Lambda_k/\lambda_k-\Omega_{k})} \Big )^{\Lambda_{k}}.
\end{equation}
   We note that $\Omega_k(\mu_N)= \omega_k$ for $1 \leq k \leq N-1$ and 
\begin{eqnarray*}
   \Omega^{\Lambda_{N}}_N(\mu_N) &=& \frac 1{ \mu^{\lambda_{N}}_N}\Big (\frac {\omega_{N-1}}{\frac {\lambda_{N}}{\Lambda_{N-1}}(\Lambda_{N-1}/\lambda_{N-1}-\omega_{N-1})} \Big )^{\Lambda_{N-1}}=\frac 1{ \mu^{\lambda_{N}}_N}\Big (\frac {\omega_{N-1}a_{N-1}}{a_{N}} \Big )^{\Lambda_{N-1}} \\
 &=& \frac 1{ \mu^{\lambda_{N}}_N}\Big (\frac {G_{N-1}}{\mu_Na_{N}} \Big )^{\Lambda_{N-1}}=\Big (\frac {G_{N}}{\mu_Na_{N}} \Big )^{\Lambda_{N}}=\Big (\frac {\Lambda_{N}}{\lambda_{N}} \Big )^{\Lambda_{N}}.
\end{eqnarray*}
   
   We now show by induction that if $\mu > e^M$, then for any $k \geq 1$,
\begin{equation}
\label{2.3}
   \Omega_k(\mu) < \frac {\Lambda_k/\lambda_k}{\Lambda_{k+1}/\lambda_{k+1}}.
\end{equation}
   As we have seen above that $\Omega_N(\mu_N)=\Lambda_{N}/\lambda_{N}$, this forces $\mu_N \leq e^M$ and hence our assertion for Theorem \ref{thm1} will follow.

   Now, to establish \eqref{2.3}, we note first the case $k=1$ follows directly from our assumption \eqref{5} on considering the case $n=1$ there. Suppose now \eqref{2.3} holds for $k \geq 1$, then by the relation \eqref{2.2}, we have
\begin{eqnarray*}
   \Omega^{\Lambda_{k+1}}_{k+1}(\mu) &=& \frac 1{ \mu^{\lambda_{k+1}}}\Big (\frac {\Omega_{k}}{\frac {\lambda_{k+1}}{\Lambda_k}(\Lambda_k/\lambda_k-\Omega_{k})} \Big )^{\Lambda_{k}} \\
&<& \frac 1{ \mu^{\lambda_{k+1}}}\Big (\frac {\frac {\Lambda_k/\lambda_k}{\Lambda_{k+1}/\lambda_{k+1}}}{\frac {\lambda_{k+1}}{\Lambda_k}(\Lambda_k/\lambda_k-\frac {\Lambda_k/\lambda_k}{\Lambda_{k+1}/\lambda_{k+1}})} \Big )^{\Lambda_{k}}=\frac 1{ \mu^{\lambda_{k+1}}}.
\end{eqnarray*}
  This implies that
\begin{equation*}
   \Omega_{k+1}(\mu) < \frac 1{ \mu^{\lambda_{k+1}/\Lambda_{k+1}}}< \frac {\Lambda_{k+1}/\lambda_{k+1}}{\Lambda_{k+2}/\lambda_{k+2}}.
\end{equation*}
  The last inequality follows from the case $n=k+1$ of our assumption \eqref{5} and this completes the proof.

\section{An Application of Theorem \ref{thm1}}
\label{sec 3} \setcounter{equation}{0}
   Our goal in this section is to establish the following:
\begin{theorem}
\label{thm2}
   Inequality \eqref{1} holds for $\lambda_k=k^{\alpha}$ for $0 < \alpha <1 $ with $C=1/(\alpha+1)$.
\end{theorem}

   We need a lemma first:
\begin{lemma}\cite[Lemma 1, 2, p.18]{L&S}
\label{lem0}
    For an integer $n \geq 1$ and $0 \leq r \leq 1$,
\begin{equation*}
    \frac {1}{r+1}n(n+1)^r \leq  \sum^n_{i=1}i^r \leq \frac {r}{r+1}\frac
   {n^r(n+1)^r}{(n+1)^r-n^r}.
\end{equation*}
\end{lemma}

  Now we return to the proof of Theorem \ref{thm2}. It suffices to check that condition \eqref{5} is satisfied with $M=1/(\alpha+1)$ there. Explicitly, we need to show that for any integer $n \geq 1$,
\begin{equation}
\label{3.1}
  \frac {\sum^n_{k=1}k^{\alpha}}{n^{\alpha}}\log \Big ( \Big ( 1+ \frac {(n+1)^{\alpha}}{\sum^n_{k=1}k^{\alpha}} \Big )\Big ( \frac {n^{\alpha}}{(n+1)^{\alpha}} \Big )\Big ) \leq \frac 1{\alpha+1}.
\end{equation}
  Now we apply Lemma \ref{lem0} to obtain:
\begin{equation*}
   1+ \frac {(n+1)^{\alpha}}{\sum^n_{k=1}k^{\alpha}}  \leq  1+ \frac {\alpha+1}{n}.
\end{equation*}
  We use this together with the upper bound in Lemma \ref{lem0} to see that inequality \eqref{3.1} is a consequence of the following inequality:
\begin{equation}
\label{3.2}
 \alpha \Big ( \log  ( 1+ \frac {\alpha+1}{n} ) - \log  (1+1/n)^{\alpha} \Big ) \leq 1-\frac 1{(1+1/n)^{\alpha}}.
\end{equation}
  We now define
\begin{equation*}
  f(x)=1-(1+x)^{-\alpha}-\alpha \Big ( \log  ( 1+ (\alpha+1)x ) -\alpha \log  (1+x) \Big ).
\end{equation*}
  Note that inequality \eqref{3.2} is equivalent to $f(1/n) \geq 0$. Hence it suffices to show that $f(x) >0$ for $0 < x \leq 1$.
Calculation shows that
\begin{equation*}
  f'(x)=\frac {\alpha g(x)}{(1+x)^{1+\alpha}\big(1+(1+\alpha)x \big )},
\end{equation*} 
  where
\begin{equation*}
  g(x)=1+ (\alpha+1)x-\big(\alpha+(1-\alpha^2)x \big )(1+x)^{\alpha}.
\end{equation*} 
   Note that when $0<\alpha <1$,
\begin{equation*}
   (1+x)^{\alpha} \leq 1+\alpha x.
\end{equation*}
    It follows that
\begin{eqnarray*}
    g(x) &\geq & 1+ (\alpha+1)x-\big(\alpha+(1-\alpha^2)x \big )(1+\alpha x) \\
         &=& 1-\alpha+\alpha x -\alpha (1-\alpha^2)x^2 :=h(x).
\end{eqnarray*}
   It is easy to see that $h(x)$ is concave for $0 \leq x \leq 1$ and $h(0)=1-\alpha >0, h(1)=1-\alpha (1-\alpha^2) >0$. It follows that $h(x) >0$ for $0<x<1$ so that $g(x) >0$ and hence $f'(x) >0$ for $0<x<1$. As $f(0)=0$, this implies $f(x) \geq 0$ for $0<x \leq 1$ and this completes the proof of Theorem \ref{thm2}.

\end{document}